\documentclass[a4paper]{article}

\usepackage{amsfonts}
\usepackage{latexsym}
\usepackage{amssymb}
\usepackage{amscd}

\newtheorem{thm}{Theorem}[section]
\newtheorem{prop}[thm]{Proposition}
\newtheorem{lemma}[thm]{Lemma}

\newenvironment{pf}{{\bf Proof.}}{\hfill$\Box$\\[1mm]}

\def\a{\alpha} \def\b{\beta} \def\g{\gamma} 
   \def\l{\lambda}
  \def\s{\sigma}  \def\t{\tau}
\def\om{\omega}

\def\bC{\mathbb{C}}
\def\bN{\mathbb{N}} 

\def\bR{\mathbb{R}}  
 
\def\bZ{\mathbb{Z}}

\def\sA{\mathcal{A}}
\def\A{\mathcal{A}}
\def\B{\mathcal{B}}
\def\Ae{\mathcal{A}^e}

\def\C{\mathcal{C}}

\def\M{\mathcal{M}}
\def\N{\mathcal{N}}

\def\sS{\mathcal{S}}

\def\V{\mathcal{V}}

\def\eacute{\mathrm{\acute{e}}}

\def\id{\mathrm{id}}

\def\im{\mathrm{im}}
\def\otherwise{\mathrm{otherwise}}
\def\mod{\mathrm{mod}}
\def\smod{\s_\mod}
\def\ker{\mathrm{ker}}

\def\Tor{\mathrm{Tor}}
\def\rank{\mathrm{rank}}

\def\Podles{\mathrm{Podle{\acute{s}}}}

\def\isom{\cong}

\def\a{\alpha}

\def\sA{{}_\s \A}
\def\AAs{\A,\sA}
\def\Asq{ \A (S^2_q )}

\def\b{\beta}

\def\suq2{SU_q (2)}
\def\sphere{S_q^2 (c,d)}
\def\tl{\triangleleft}

\def\Asn{{(\A,\s)}^{\natural}}
\def\Asn{\A^{\natural}_\s}
\def\Asnzero{\A^{\natural}_{\s,0}}
\def\Asnone{\A^{\natural}_{\s,1}}
\def\Asntwo{\A^{\natural}_{\s,2}}
\def\Asnthree{\A^{\natural}_{\s,3}}
\def\Asnn{\A^{\natural}_{\s,n}}

\def\Asn{{(\A,\s)}^{\natural}}
\def\Asn{C^\s}
\def\Asnzero{C^\s_0}
\def\Asnone{C^\s_1}
\def\Asntwo{C^\s_2}
\def\Asnthree{C^\s_3}
\def\Asnn{C^\s_n}

\def\Acd{\A(c,d)}

\def\csuq2{C(SU_q (2))}
\def\slq2{SL_q (2)}
\def\s{\sigma}

\def\rto{\rightarrow}

\def\z{\zeta}

\def\blfootnote{\xdef\@thefnmark{}\@footnotetext}

\begin{document}

\title{Twisted cyclic homology of all $\Podles$ quantum spheres}
\author{Tom~Hadfield\footnote{Supported until 31/12/2003 by the EU Quantum Spaces - Noncommutative Geometry Network (INP-RTN-002) and from 1/1/2004 by an EPSRC postdoctoral fellowship}}
\date{\today}
\maketitle

\centerline{School of Mathematical Sciences}
\centerline{Queen Mary, University of London}
\centerline{327 Mile End Road, London E1 4NS, England}
\centerline{t.hadfield@qmul.ac.uk}
\centerline{MSC 2000; 58B34, 19D55, 81R50, 46L} 
\centerline{Keywords: cyclic homology, Hochschild homology,}
\centerline{quantum group, quantum sphere}

\begin{abstract}
We calculate the twisted Hochschild and cyclic homology of all Podle\'s quantum spheres relative to arbitrary  automorphisms.
  The dimension drop in Hochschild homology is overcome via twisting by  the modular automorphism of the 
  canonical $SU_q (2)$-invariant linear functional.
Specializing  to the standard quantum sphere, we identify the cohomology class of the 2-cocycle discovered by Schm\"{u}dgen and Wagner corresponding to the distinguished covariant differential calculus found by Podle\'s.
\end{abstract}

\maketitle

\section{Introduction}

Twisted cyclic cohomology was discovered by Kustermans, Murphy and Tuset \cite{kmt}, arising naturally from covariant differential calculi over compact quantum groups. 
They defined a cohomology theory relative to a pair of  an algebra $\A$ and an automorphism $\s$, which on taking $\s = \id$ reduces to ordinary cyclic cohomology of $\A$. 
While it was immediately recognised that twisted cyclic cohomology (and its dual, twisted cyclic homology, the subject of this paper) fits into Connes' general framework of cyclic objects, its relation with 
 differential calculi \cite{sw1,sw2} and recent connection with the ``dimension drop" phenomenom in Hochschild homology \cite{ft,hk,hk2,si} makes it of independent interest .

Previously \cite{hk} we studied the twisted Hochschild and cyclic homology of the quantum $SL(2)$ group . 
We now extend this work to the $\Podles$ quantum spheres \cite{pod87,pod92}, which are ``quantum homogeneous spaces'' for quantum $SL(2)$. 
The $\Podles$ spheres have been extensively studied, with much work done  
 constructing Dirac operators, spectral triples and the corresponding local index formulae. We mention only \cite{bk,cpb,ds,nt} amongst many others. 
In general, covariant differential calculi over quantum groups do not fit into Connes' formalism of spectral triples \cite{schm}.
However, in \cite{sw2} Schm\"udgen and Wagner constructed a Dirac operator giving a commutator representation of the distinguished 2-dimensional first order covariant calculus over the $\Podles$  sphere \cite{pod92}. The associated twisted cyclic 2-cocycle $\tau$ was shown to be a nontrivial element of twisted cyclic cohomology. This 2-cocycle does not correspond to the ``no dimension drop" case - the fact that twisting overcomes the dimension drop in Hochschild homology for the Podle\'s spheres is the main new result of this paper.

A summary of this paper is as follows.
In section \ref{section:twisted_hoch_and_cyclic} we recall the definitions  \cite{hk}, \cite{kmt} of twisted Hochschild and cyclic homology $HH_\ast^\s (\A)$, $HC_\ast^\s (\A)$.
These ``twisted homologies" arise from a cyclic object in the sense of Connes \cite{loday}, hence all Connes' homological machinery can be applied.
Previously we proved that:
\begin{thm}
\label{crucial_isomorphism}
 \cite{hk}
For arbitrary $\A$ and $\s$, if $\s$ acts diagonally relative to a set of generators of $\A$ then
 $HH_n^\s (\A) \isom H_n (\AAs)$ for each $n$.
\end{thm}
Here $\sA$ is the  ``$\s$-twisted'' $\A$-bimodule with $\A$ as  underlying vector space, and $\A$-bimodule structure 
\begin{equation}
\label{twisted_bimodule}
a_1  \triangleright x \triangleleft a_2 = \s(a_1) x a_2 \quad  x, a_1, a_2 \in \A
\end{equation}
It was shown by Kr\"ahmer \cite{ulrich_thesis} that all automorphisms of the $\Podles$  spheres are diagonalisable, hence using 
 $H_n (\AAs)$ $\isom \Tor^{\Ae}_n (\sA,\A)$ \cite{loday} ($\A^e = \A \otimes \A^{op}$), if we have a projective resolution of $\A$ by left $\A^e$-modules, we can in principle compute $HH_n^\s (\A)$.

Hochschild and cyclic homology of the $\Podles$ quantum spheres
was calculated by Masuda, Nakagami and Watanabe  \cite{mnw}, 
using a free resolution that we rely on in this paper. 
In section \ref{section:podles} we recall their definitions. 
  In section \ref{section_H_*(A,A,s)} we use this resolution to calculate the 
Hochschild homologies $H_n (\AAs)$, which by Theorem \ref{crucial_isomorphism} are isomorphic to the twisted Hochschild homologies $HH_n^\s (\A)$.

We obtain the following striking result (Theorem \ref{thm_H_2_nonzero}).
In the untwisted situation \cite{mnw} the Hochschild groups $HH_n (\A) = H_n (\A,\A)$ vanish for $n \geq 2$, in contrast to the classical situation $q=1$ (the ordinary 2-sphere) whose Hochschild dimension is 2. 
This ``dimension drop'' phenomenon has been seen in many other quantum situations (see \cite{ft} for a detailed discussion).
However, in the twisted situation, there exist automorphisms $\s$ with $HH_n^\s (\A) \neq 0$ for $n=0,1,2$. 
These automorphisms are  precisely the positive powers of the canonical modular automorphism associated to the $SU_q (2)$-invariant linear functional discovered by Noumi and Mimachi \cite{mn}. 
 For the standard quantum sphere, which naturally embeds as a subalgebra of quantum $SU(2)$, this modular automorphism coincides with the modular automorphism induced from the Haar state on quantum $SU(2)$. 
 The central role of the modular automorphism in avoiding the dimension drop in Hochschild homology was also  seen for quantum $SL(N)$ \cite{hk,hk2}.
 Similar results have been obtained by Sitarz \cite{si} for quantum hyperplanes.

In section \ref{section:twisted_cyclic_homology} we calculate twisted cyclic homology as the total homology of Connes' mixed $(b,B)$-bicomplex arising from the underlying cyclic object.
Finally, in section \ref{standard_quantum_sphere} we apply our results to the standard quantum sphere, showing that  the class $[\t]$ in twisted cyclic cohomology $HC^2_\s (\A)$ of  Schm\"{u}dgen and Wagner's  twisted cyclic 2-cocycle is proportional to  $[ Sh_A]$, where $S$ is the periodicity operator and $h_A$ an explicit  nontrivial twisted cyclic 0-cocycle.

\section{Twisted Hochschild and cyclic homology}
\label{section:twisted_hoch_and_cyclic}

We recall the definitions of twisted Hochschild and  cyclic homology \cite{hk}.
Let $\A$ be a unital algebra over a field $k$ (assumed to be of characteristic zero), and $\s$ an automorphism.
Define $C_n (\A) = \A^{\otimes(n+1)}$. For brevity, we will write
$a_0 \otimes  \ldots \otimes a_n  \in \A^{\otimes (n+1)}$ as $(a_0, \ldots , a_n )$.
Define the twisted cyclic operator $\l_\s : C_n (\A) \rto C_n (\A)$ by
$\l_\s ( a_0, \ldots  , a_n ) = (-1)^n (\s(a_n), a_0,  \ldots ,a_{n-1})$.
Hence
$\l^{n+1}_\s  ( a_0, \ldots  , a_n ) = ( \s(a_0 ) , \ldots , \s( a_n))$.
 Now consider the quotient
 \begin{equation}
 \label{C_n^sigma}
 C_n^\s (\A) = \A^{\otimes(n+1)} / (\id - {\l_\s^{n+1}})
 \end{equation}
 If $\s = \id$, then $C_n^\s (\A) = \A^{\otimes(n+1)}$.
The twisted Hochschild boundary operator
$b_\s : C_{n+1} (\A) \rightarrow C_n (\A)$  is given by 
$$ b_\s ( a_0, \ldots , a_{n+1} ) =
\sum_{j=0}^{n} {(-1)}^j (a_0, ., a_j a_{j+1}, .,a_{n+1} )$$
\begin{equation}
\label{b_sigma}
+ {(-1)}^{n+1} ( \s(a_{n+1}) a_0, a_1 , \ldots, a_n)
\end{equation}
We have $b_\s^2 = 0$ and $b_\s \; \l_\s^{n+2} = \l_\s^{n+1} \; b_\s$, hence $b_\s$ descends to the quotient,\\
 $b_\s : C_{n+1}^\s (\A) \rightarrow C_n^\s (\A)$.
Twisted Hochschild homology $HH_*^\s (\A)$ is defined as the homology of the complex
 ${\{ C_n^\s (\A), b_\s \}}_{n \geq 0}$. Taking $\s = \id$ gives $HH_* (\A) = H_* (\A,\A )$, the Hochschild homology of $\A$ with coefficients in $\A$.

Now define $C_n^{\s,\l} (\A) = \A^{\otimes(n+1)} / (\id - \l_\s)$. We have a surjection
$C_n^\s (\A) \rto C_n^{\s,\l} (\A)$.
As maps $\A^{\otimes (n+1)} \rto \A^{\otimes n}$, we have
$b_\s ( \id - \l_\s) = (\id - \l_\s) b'$, where
\begin{equation}
\label{b_prime}
b' ( a_0, \ldots , a_n) =\sum_{j=0}^{n} {(-1)}^j (a_0, ,\dots , a_j a_{j+1}, \ldots ,a_n )
\end{equation}
Hence $b_\s$ descends to a map
$b_\s : C_{n+1}^{\s,\l} (\A) \rto C_n^{\s,\l} (\A)$.
Twisted cyclic homology  $HC_*^\s (\A)$ is then defined as the homology of the complex
 ${\{ C_n^{\s,\l} (\A), b_\s  \}}_{n \geq 0}$.
Taking $\s = \id$ gives back ordinary cyclic homology $HC_* (\A)$.

Equivalently,  twisted cyclic homology is the total homology of Connes' mixed $(b,B)$-bicomplex coming from the underlying cyclic object, which we define in section \ref{section:twisted_cyclic_homology}, and use to calculate $HC_*^\s (\A)$ from $HH_*^\s (\A)$ for the $\Podles$ spheres.

\section{The Podle\'s quantum spheres}
\label{section:podles}

\subsection{\sc the coordinate algebras $\Acd$}

Let $k$ be a field of characteristic zero, and $q \in k$ nonzero and not a root of unity. 
For  $c$, $d \in k$, with $c + d \neq 0$,  we define the coordinate algebra $\A ( c,d)$ of the $\Podles$ quantum two sphere $\sphere$ to be the unital  $k$-algebra with generators  $A$, $B$, $B^\ast$ satisfying 
\begin{equation}
\label{BA=}
BA = q^2 AB, \quad A B^\ast = q^2 B^\ast A
\end{equation}
$$B^* B = cd + (c-d) A - A^2, \quad
B B^*  = cd + q^2 (c-d) A - q^4 A^2$$
In the notation of \cite{mnw}, we take $A = \z$, $B = Y$, $B^* = -qX$. As algebras, $\A( rc, rd) \isom \A(c,d)$ for any  $r \in k$, $r \neq 0$.
A Poincar${\eacute}$-Birkhoff-Witt basis for $\Acd$ consists of the monomials
\begin{equation}
\label{pbw_basis}
{\{ B^j A^k  \}}_{j,k \geq 0}, \quad {\{  {B^\ast}^{j+1} A^k \}}_{j,k \geq 0}
\end{equation}
Working over $\bC$ (we take $q, c, d \in \bR$, with $0 < q<1$, $0<c$), for $d >0$, there is a family of quantum spheres parameterised by $t \in \bR$, $t >0$, with
$$B^* B = t1 + A - A^2,\quad B B^* = t1 + q^2 A - q^4 A^2$$
and also the ``equatorial quantum sphere'', with $B^* B =  1 - A^2$, $B B^* = 1 - q^4 A^2$.
The C*-algebraic completions (with $A^\ast = A$) 
were shown by Sheu \cite{sheu} to all be isomorphic.
However,  Kr\"ahmer proved  the underlying algebras are pairwise non-isomorphic \cite{ulrich_thesis}.
Taking $t=0$ gives the ``standard quantum 2-sphere"
\begin{equation}
\label{standard_qsphere}
B^* B =  A - A^2,\quad B B^* =  q^2 A - q^4 A^2
\end{equation}
Now recall that the coordinate Hopf *-algebra $\A(SU_q (2))$ is the unital *-algebra over $\bC$ (algebraically) generated by elements $a$, $c$ satisfying the relations 
$$
a^* a + c^* c  =1,\quad
a a^* + q^2 c^* c  = 1, \quad c^* c = c c^*,\quad
ac = qca,\quad
a c^* = q c^* a
$$
There is a dual pairing $< .,.>$ of $\A(SU_q (2))$ with  $U_q ( su(2))$, with standard generators $E$, $F$, $K^{\pm 1}$ \cite{sw2}, giving
 left and  right actions of $U_q ( su(2))$ on $\A(SU_q (2))$:
\begin{equation}
\label{right_action}
\label{left_action}
f \triangleright x = \sum \; < f, x_{(2)} > \; x_{(1)}, \quad
x \triangleleft f = \sum \; <f , x_{(1)} > \; x_{(2)}
\end{equation}
The coordinate *-algebra $\Asq$ of the standard Podle\'s quantum sphere is the *-subalgebra of $\A(SU_q (2))$ invariant under the action of the grouplike element 
$K \in U_q ( su(2))$.
Explicitly,
$$a \triangleleft K = q^{-1/2} a, \quad
a^* \triangleleft K = q^{1/2} a^*, \quad
c \triangleleft K = q^{1/2} c, \quad
c^* \triangleleft K = q^{-1/2} c^*$$
Writing $A = c^* c$, $B = ac$, $B^* = c^* a^*$ gives  the relations (\ref{BA=}), (\ref{standard_qsphere}).

Masuda, Nakagami and Watanabe \cite{mnw} gave a resolution of $\A = \Acd$,
\begin{equation}
\label{mnw_resolution}
\ldots  \rto \M_{n+1} \rto \M_n \rto \ldots \rto \M_2 \rto \M_1 \rto \M_0 \rto \A \rto 0
\end{equation}
by free left $\A^e$-modules $\M_n$ ($\A^e = \A \otimes \A^{op}$), with 
$\rank(\M_0) =1$, $\rank(\M_1) =3$, $\rank(\M_n) = 4$ for $n \geq 2$.
Adapting their notation, 
$\M_1$ has a basis $\{ e_A, e_B, e_{B^\ast} \}$, 
with $d_1 :\M_1 \rto \M_0 = \A^e$ given by 
\begin{equation}
\label{defn_d_1}
d_1 ( e_t) = t \otimes 1 - 1 \otimes t^o, \quad t = A, B, B^\ast 
\end{equation}
$\M_2$ has basis $\{ e_A \wedge e_B$, $e_A \wedge e_{B^\ast}$, 
$\vartheta_S^{(1)}$,
$\vartheta_T^{(1)} 
\}$, with $d_2 : \M_2 \rto \M_1$ given by
$$d_2 ( 1_{\A^e}  \otimes (e_A \wedge e_{B^\ast}) ) = 
(A \otimes 1 - 1 \otimes q^2 A^o ) \otimes e_{B^\ast} - (q^2 B^\ast \otimes 1 - 1 \otimes {B^\ast}^o) \otimes e_A$$
$d_2 ( 1_{\A^e} \otimes (e_A \wedge e_B) ) =  ( q^2 A \otimes 1 - 1 \otimes A^o ) \otimes e_B - (B \otimes 1 - 1 \otimes q^2 B^o ) \otimes e_A$
$$d_2 ( 1_{\A^e} \otimes \vartheta_S^{(1)}) =  -q^{-1} \{ B \otimes 1 \otimes e_{B^\ast} + 1  \otimes {B^\ast}^o \otimes e_B \}
- q \{ q^2 ( A \otimes 1 + 1 \otimes A^o ) - (c-d) \} \otimes e_A$$
\begin{equation}
\label{defn_d2}
d_2 ( 1_{\A^e} \otimes \vartheta_T^{(1)} ) = -q^{-1} \{ 1 \otimes B^o \otimes e_{B^\ast} + B^\ast \otimes 1 \otimes e_B \} - q^{-1} \{ (A \otimes 1 + 1 \otimes A^o) - (c-d) \} \otimes e_A \nonumber
\end{equation}
$\M_3$ has basis $\{ e_A \wedge \vartheta_S^{(1)}$, 
 $e_A \wedge \vartheta_T^{(1)}$, $e_{B^\ast} \wedge \vartheta_S^{(1)}$,
 $e_B \wedge \vartheta_T^{(1)} \}$, with $d_3 : \M_3 \rto \M_2$ 
 $$d_3 ( 1_{\A^e} \otimes (e_A \wedge \vartheta_S^{(1)})) = 
 (A \otimes 1 - 1 \otimes A^o) \otimes \vartheta_S^{(1)} + q^{-3} \{ B \otimes 1 \otimes (e_A \wedge e_{B^\ast}) + 1 \otimes {B^\ast}^o \otimes (e_A \wedge e_B) \}$$
 $$d_3 (1_{\A^e} \otimes (e_A \wedge \vartheta_T^{(1)}))= 
 (A \otimes 1 - 1 \otimes A^o) \otimes \vartheta_T^{(1)} + q^{-1} \{ 1 \otimes B^o \otimes (e_A \wedge e_{B^\ast}) + B^\ast \otimes 1 \otimes (e_A \wedge e_B)\}$$\\
 $d_3 (1_{\A^e} \otimes (e_{B^\ast} \wedge \vartheta_S^{(1)}))= 
B^\ast \otimes 1 \otimes \vartheta_S^{(1)} - 1 \otimes {B^\ast}^o \otimes \vartheta_T^{(1)}$\\
$$ -q^{-1} \{ ( A \otimes 1 + 1 \otimes q^2 A^o ) - (c-d) \} \otimes (e_A \wedge e_{B^\ast} )$$\\
$d_3( 1_{\A^e} \otimes (e_B \wedge \vartheta_T^{(1)}))= B \otimes 1 \otimes \vartheta_T^{(1)} - 1 \otimes B^o \otimes \vartheta_S^{(1)}$
\begin{equation}
\label{defn_d3}
 - q^{-1} \{ ( q^2 A \otimes 1 + 1 \otimes A^o ) - (c-d)  \} \otimes (e_A \wedge e_B)
 \end{equation}
$\M_4$ has basis $\{ e_A \wedge e_{B^\ast} \wedge \vartheta_S^{(1)}$,  $e_A \wedge e_B \wedge \vartheta_T^{(1)}$, $\vartheta_S^{(2)}$, $\vartheta_T^{(2)} \}$, with $d_4 : \M_4 \rto \M_3$\\
$d_4 ( 1_{\A^e} \otimes (e_A \wedge e_{B^\ast} \wedge \vartheta_S^{(1)} )) =$
$$( A \otimes 1 - 1 \otimes q^2 A^o ) \otimes ( e_{B^\ast} \wedge \vartheta_S^{(1)})
- q^2 B^\ast \otimes 1 \otimes ( e_A \wedge \vartheta_S^{(1)}) + 1 \otimes {B^\ast}^o \otimes ( e_A \wedge \vartheta_T^{(1)})$$
$d_4 ( 1_{\A^e} \otimes ( e_A \wedge e_B \wedge \vartheta_T^{(1)}) =$
$$( q^2 A \otimes 1 - 1 \otimes A^o ) \otimes ( e_B \wedge \vartheta_T^{(1)})
- B \otimes 1 \otimes ( e_A \wedge \vartheta_T^{(1)}) + 1 \otimes q^2 B^o \otimes ( e_A \wedge \vartheta_S^{(1)})$$
$d_4 ( 1_{\A^e} \otimes  \vartheta_S^{(2)} ) = - q^{-1} B \otimes 1 \otimes ( e_{B^\ast} \wedge \vartheta_S^{(1)})) - q^{-1} \otimes {B^\ast}^o \otimes ( e_B \wedge \vartheta_T^{(1)})$
$$- q [ q^2 ( A \otimes 1 + 1 \otimes A^o ) - (c-d) ] \otimes ( e_A \wedge \vartheta_S^{(1)}))$$
$d_4 ( 1_{\A^e} \otimes \vartheta_T^{(2)} ) = - q^{-1} \otimes B^o \otimes ( e_{B^\ast} \wedge \vartheta_S^{(1)}) - q^{-1} B^\ast \otimes 1 \otimes ( e_B \wedge  \vartheta_T^{(1)}) $
$$- q^{-1} ( (A \otimes 1 + 1 \otimes A^o ) - (c-d) ) \otimes ( e_A \wedge \vartheta_T^{(1)})$$

  We refer the reader to \cite{mnw} for the  $\M_n$ and $d_n$ for $n \geq 5$.
  In section \ref{section_H_*(A,A,s)} we use this resolution to 
 calculate the Hochschild homology 
 $H_* (\AAs)$ of $\A =\Acd$ with coefficients in the twisted $\A$-bimodule ${}_\s \A$ defined in (\ref{twisted_bimodule}).

\subsection{\sc comparison of the m-n-w and bar resolutions}

We wish to identify generators of $H_* (\AAs)$, found as elements of the modules $\M_n$, with Hochschild cycles realised as elements of $\A^{\otimes n}$. 
Recall \cite{loday} the bar resolution, with  differential $b'$ given by (\ref{b_prime})
$$\ldots \rto \A^{\otimes (n+2)} \rto^{b'}
\A^{\otimes (n+1)} \rto \ldots \rto
\A^{\otimes 2} \rto^{b'} \A \rto 0$$
which is a projective resolution of $\A$ as a left $\A^e$-module.
Each $\A^{\otimes (n+1)}$ is a left $\A^e$-module via 
$(x \otimes y^o )(a_0, a_1, \ldots a_n ) =
(x a_0, a_1 , \ldots , a_n y )$.
The comparison theorem (see, for example \cite{w}, Theorem 2.2.6) says that given a projective resolution $\ldots \rto \M_1 \rto^{d_1} \M_0 \rto^{d_0} \B \rto 0$ of a left $\A$-module $\B$, and a map $f:\B \rto \C$, then for every resolution $\ldots \rto \N_1 \rto \N_0 \rto^\eta \C \rto 0$ there is a chain map $\{ f_i : \M_i \rto \N_i \}_{i \geq 0}$, unique up to chain homotopy equivalence, lifting $f$ in the sense that $\eta \circ f_0 = f \circ d_0$.
In our situation, taking $\B=\C=\A$ and $f=\id$, maps $f_i : \M_i \rto \A^{\otimes (i+2)}$ giving a 
 commutative diagram

	\[
     	\begin{CD}
   	{\ldots} @>>> {\M_3}  @> {d_3} >> {\M_2} @ >{d_2} >> {\M_1} @ >{d_1} >> {\M_0} @ >{d_0} >> {\A} @ >>> {0}\\
	@ . @ V{f_3} VV @ V{f_2} VV @ V{f_1} VV @ V{f_0} VV @ V{\isom} VV @ .\\
	{\ldots} @>>> {\A^{\otimes 5}} @> {b'} >> {\A^{\otimes 4}} @ >{b'}>> {\A^{\otimes 3}} @ >{b'} >> {\A^{\otimes 2}} @ >{b'} >> {\A} @ >>> {0}\\
     	\end{CD}
	\]\\
 are given by, in the notation of the previous section:
\begin{eqnarray}
\label{defn_f_0}
\label{defn_f_1}
\label{defn_f_2}
&& f_0 (a_1 \otimes a_2^o ) = (a_1 , a_2), \quad
 f_1 ( e_t ) = (1,t,1), \quad t = A, \, B,\, B^\ast \nonumber\\
&& f_2 ( e_A \wedge e_{B^\ast} ) =
(1, A, B^\ast ,1) - q^2 (1,B^\ast ,A,1) \nonumber\\
&& f_2 ( e_A \wedge e_B ) = 
q^2 (1,A,B,1) - (1,B,A,1) \nonumber\\
&& f_2 ( \vartheta^{(1)}_S ) = - q^{-1} (1,B,B^\ast ,1) - q^3 (1,A,A,1) - q^{-1} cd (1,1,1,1) \nonumber\\
&& f_2 ( \vartheta^{(1)}_T ) = -q^{-1} (1,B^\ast ,B,1) - q^{-1} (1,A,A,1) -q^{-1} cd (1,1,1,1)
\end{eqnarray}
Higher $f_i$ can be found inductively: the above is as much as we will need in the sequel.
Applying ${}_\s \A \otimes_{\A^e} -$  to both resolutions allows us to identify generators of homology found from the M-N-W resolution with explicit Hochschild cycles.

\subsection{\sc automorphisms of $\Acd$}

It was shown by Kr\"ahmer \cite{ulrich_thesis} that every automorphism of $\A(c,d)$ acts 
 diagonally with respect to the generators $A$, $B$, $B^\ast$.
For $c \neq d$, every  automorphism  is of the form 
\begin{equation}
\label{defn_sigma}
\s_\l (B) = \l B, \quad \s_\l (A) = A, \quad \s_\l ( B^\ast ) = \l^{-1} B^\ast
\end{equation}
some $\l \in k$, $\l \neq 0$.  If $c=d$, there is a second family of automorphisms
\begin{equation}
\label{defn_tau_aut}
\t_\l (B) = \l B, \quad \t_\l (A) = -A, \quad \t_\l ( B^\ast ) = \l^{-1} B^\ast
\end{equation}
It follows from Theorem \ref{crucial_isomorphism} that:

\begin{lemma}
\label{cor_crucial_isomorphism}
$HH_n^\s (\A) \isom H_n (\A, \sA)$ for all $n \geq 0$ and every $\s$.
\end{lemma}

Working over $\bC$, Noumi and Mimachi \cite{mn} proved the existence of a unique linear functional $h : \A(c,d) \rto \bC$ invariant under the left coaction of quantum $SU(2)$, and satisfying $h(1)=1$. On monomials this is given by 
\begin{equation}\label{nm_functional}
h( B^{m+1} A^n) = 0 = h( (B^\ast )^{m+1} A^n ), \quad h( A^n ) = {\tiny{\frac{f(0)}{f(n)}}} ( \;{\frac{c^{n+1} - (-d)^{n+1}}{c+d}} \;)
\end{equation}
where $f(n) = q^{-2} - q^{2n}$.
$h$ is a twisted cyclic 0-cocycle. 
Borrowing terminology used for  quantum $SU(2)$,
 the unique automorphism $\smod$ satisfying $h(xy) = h(y \smod (x))$ is called 
 the modular automorphism (so $h$ is a $\smod^{-1}$-twisted 0-cocycle). Concretely, 
\begin{equation}
\label{mod_aut_podles}
\smod : \quad A \mapsto A, \quad B \mapsto q^{-2} B, \quad B^\ast \mapsto q^2 B^\ast
\end{equation}
Obviously $\smod$  is well-defined over any field, not just $\bC$. As previously seen, 
the standard Podle\'s quantum sphere  embeds as a subalgebra of quantum $SU(2)$, and the modular automorphism associated to the  Haar state on quantum $SU(2)$ restricts to an automorphism of the standard Podle\'s sphere coinciding with (\ref{mod_aut_podles}).

\section{Twisted Hochschild homology} 
\label{section_H_*(A,A,s)}

We calculate the Hochschild homologies $H_n (\A, \sA)$ of $\A = \Acd$ for all  automorphisms $\s = \s_\l$, $\t_\l$ using the Masuda-Nakagami-Watanabe resolution (\ref{mnw_resolution}). 
By Lemma \ref{cor_crucial_isomorphism} we can identify these with $HH_n^\s (\A)$. 
The case $\s = \id$ was already treated in \cite{mnw}. 
In each case we exhibit explicit generators.

\subsection{ $HH_0^\s (\A)$}

Let $\s_\l$, $\t_\l$ be the automorphisms of $\Acd$ given by (\ref{defn_sigma}), (\ref{defn_tau_aut}).

\begin{prop}
\label{HH_0}
For arbitrary $c$ and $d$ (with $c+d \neq 0$) and $\s = \s_\l$ we have:
\begin{enumerate}
\item{ For $\l=1$ ($\s =\id$), $HH_0^\s (\A)$ is countably infinite dimensional.}
\item{For $\l \neq 1$, 
$HH_0^\s (\A) \isom k^2$.}
\end{enumerate}
For $c=d$, and $\s = \t_\l$ we have:
\begin{enumerate}
\item{ For $\l=1$, $HH_0^\s (\A)$ is countably infinite dimensional.}
\item{For $\l \neq 1$, 
$HH_0^\s (\A) \isom k$.}
\end{enumerate}
  \end{prop}

\begin{pf}
We have 
$HH_0^\s (\A) = \{ \; [a] \; : \; a = \s(a), \; [a_1 a_2 ] = [\s(a_2) a_1] \; \}$.
Hence for $\s = \s_\l$  with $\l \neq 1$, we need only consider P-B-W monomials $A^n$.
Now,
$$cd [ A^n] = [ (A^2 + (d-c) A + B^\ast B) A^n ]
= [ A^{n+2} ] + (d-c) [ A^{n+1}] + q^{2n} [  \s(B) B^\ast A^n ],$$
$$\Rightarrow [ A^{n+2}] + (d-c) [A^{n+1}]  -cd [A^n ] 
= q^{2n}  \l ( q^4 [A^{n+2}] + (d-c) q^2 [A^{n+1}] -cd  [A^n]),$$
so $f(n+2) [A^{n+2}] + (d-c) f(n+1) [A^{n+1}] -cd f(n) [ A^n] =0$, where $f(n) = \l^{-1} - q^{2n}$.
Write $x_n = f(n) [A^n]$. Then we have
\begin{equation}
\label{x_n_recurrence_relation}
x_{n+2} + (d-c) x_{n+1} -cd x_n =0 \quad \; \forall \; n \geq 0
\end{equation}
For $\l \notin q^{-2 \bN}$, we have
$x_n = (c+d)^{-1} (\a c^n + \b (-d)^n )$
with $\a$, $\b$ given by:
$$\a = d f(0) [1] + f(1) [A], \quad \b = c f(0) [1] - f(1) [A]$$
So for $\l \notin q^{-2 \bN}$, we have 
$HH_0^\s ( \A) \isom k [1] \oplus k [A]$. There are three remaining cases we treat seperately:\\

{\bf Case 1:} $\s = \s_\l$, $\l = q^{-(2b+2)}$ ($b \geq 0$). 
Solving (\ref{x_n_recurrence_relation}) 
 requires some care.
However, it is not difficult to show that:
\begin{enumerate}
\item $c \neq d$. $HH_0^\s (\A) \isom k[1] \oplus k[ A^{b+1}]$. If $c \neq d$ then $[A]$, $[A^{b+1}]$ also span.
\item $c=d$. If $\l = q^{-(4b+2)}$, then  
$HH_0^\s (\A) \isom k[1] \oplus k[A^{2b+1}]$.\\
For $\l = q^{-(4b+4)}$, $HH_0^\s (\A) \isom k[A] \oplus k[ A^{2b+2} ]$.
\end{enumerate}

We give the proof of case 1 (case 2 is  similar). 
For $\l = q^{-(2b+2)}$, $f(b+1) =0$, hence $x_{b+1} =0$. So (\ref{x_n_recurrence_relation}) holds for $n \neq b$, $b \pm 1$. Hence for $n \geq b+2$ we have
$x_n = (c+d)^{-1} ( \a c^n + \b (-d)^n )$ with 
$$\a =  c^{-(b+2)} [ x_{b+3} + d x_{b+2} ], \quad \b =  (-d)^{-(b+2)}[ c x_{b+2} - x_{b+3}]$$
Further, we have $x_{b+3} + (d-c) x_{b+2} =0$, $x_{b+2} -cd x_b =0$, 
so $x_{b+2} = cd x_b$, $x_{b+3} = cd(c-d) x_b$, hence
$\a = c^{-b} d x_b$, $\b = c (-d)^{-b} x_b$. Also, for $b \geq 1$ we have $(d-c) x_b - cd x_{b-1} =0$.
Finally, for $0 \leq n \leq b-2$ (provided $b \geq 2$) (\ref{x_n_recurrence_relation}) holds, and solving this gives $x_n$ for each $n \leq b$ in terms of $x_b$. 
We have, for each $b \geq 0$,
$$x_n = g(n-b-1) x_b,  \quad \forall \; n \geq 0$$
where for  $t \in \bZ$, $g(t) = (c+d)^{-1} cd [ c^t - (-d)^t ]$.
So for $cd =0$, $[A^n ] =0$ for $n \neq 0$, $b+1$, while for $cd \neq 0$
each $x_n$, for $n \neq b+1$, is a nonzero multiple of $x_b$, and so of $x_0$. Since $f(n) \neq 0$ for $n \neq b+1$, we have $[ A^n ] = \rho_n [1]$, some $\rho_n \neq 0$,  for each $n \neq b+1$.
So for  $b \geq 0$, $[1]$, $[A^{b+1}]$, equivalently (for $b \geq 1$) $[A]$, $[A^{b+1}]$, span $HH_0^\s (\A)$. 
For nontriviality and linear independence, consider $\s$-twisted 0-cocycles $\t_0$, $\t_{b+1}$, defined (for $cd \neq 0$) on monomials $x = B^m A^n$ by
$$\t_0 (x) =\left\{
	\begin{array}{ll}
	{\tiny {\frac{g(n-b-1)}{f(n)}}} \; &: \; x = A^n, n \neq b+1\\
	0 \; &: \; \otherwise
	\end{array}\right\},\nonumber\\
	\t_{b+1} (x) =\left\{
	\begin{array}{ll}
	1 \; &: \; x = A^{b+1}\\
	0 \; &: \; \otherwise
	\end{array}\right.\nonumber\\$$
	For $cd =0$, define $\t_0 (1) =1$, $\t_0 (x) =0$ otherwise. 
Then for all $c \neq d$, $\t_0 (1)  \neq 0$, $\t_0 (A^{b+1}) = 0$.
So $HH_0^\s (\A) = k^2$, with  basis $[1]$, $[A^{b+1}]$.
We note the similarity of $\t_0$ with Noumi and Mimachi's  $SU_q (2)$-invariant functional $h$ (\ref{nm_functional}),   although the latter corresponds to the case $\l =q^2$.\\ 	

 {\bf Case 2:} $\s = \s_\l$, $\l = 1$ ($\s = \id$). We have $x_0 =0$, and: 
\begin{enumerate}
\item{$cd =0$, $c-d \neq 0$ : $x_{n+1} = (c-d)^n x_1$ for all $n \geq 0$.}
\item{$cd \neq 0$, $c-d =0$ : $x_{2n+1} = (cd)^n x_1$, $x_{2n+2}=0$, for all $n \geq 0$.
}
\item{$cd \neq 0$, $c-d \neq 0$ : Then $x_{n+1} = g(n) x_1$, for some function $g$.}
\end{enumerate}
Also 
$[ A^m B^n ] = [ \s( B^s) A^m B^{n-s} ] = q^{2sm} [ A^m B^n ]$ for $0 \leq s \leq n$.
So $[A^m B^n ] =0$ unless $m=0$ or $n=0$. 
Similarly for $[ A^m {B^\ast}^n]$. 
So for $\s = \id$, exactly as in \cite{mnw},
\begin{equation}
\label{HH_0_id}
HH_0^{\id} (\A) = H_0 ( \A,\A) \isom k [1] \oplus k [A] \oplus \; ( \Sigma_{m \geq 1}^\oplus \; k [B^{m}] ) \; \oplus \; ( \Sigma_{m \geq 1}^\oplus \; k [{B^\ast}^{m}] )
\end{equation}

{\bf Case 3:} $c=d$, $\s = \t_\l$. Then 
$[A^{n+1} ] = [ A^n A] = [ \s(A) A^n ] = - [A^{n+1} ]$. So $[A^{n+1} ] =0$ for $n \geq 0$.
So for $\l \neq 1$, $HH_0^\s (\A) \isom k[1]$, and for $\l=1$,
$HH_0^\s (\A)$ is given by (\ref{HH_0_id}), except that $[A]=0$.
\end{pf}

\subsection{ $HH_1^\s (\A)$}

\begin{prop} \label{HH_1}
For $\s = \t_\l$,  if $\l \neq 1$ then $HH_1^\s (\A) =0$.
 For $\l = 1$, $HH_1^\s (\A)$ is countably infinite dimensional, spanned by
 $[( B^j , B)]$, $[( {B^\ast}^j , B^\ast )]$,  $j \geq 0$.\\
 For $\s = \s_\l$, and arbitrary $c$ and $d$ (with $c+d \neq 0$) we have
\begin{enumerate}
\item For  $\l = q^{-2}$ or $\l \notin q^{-2 \bN}$,  $HH_1^\s (\A)  \isom k [ (1,A)]$.
\item For  $\l = 1$ ($\s = \id$), $HH_1^\s (\A)$ is countably infinite dimensional, spanned by
$[(1,A)]$, $[(B^j , B)]$, $[( {B^\ast}^j , B^\ast )]$  ($j \geq 0$).
\item For  $cd =0$, and $\l = q^{-(2b+4)}$ ($b \geq 0$), $HH_1^\s (\A)  \isom k [ (A^{b+1}, A)]$.
\item For  $c-d=0$, if $\l = q^{-(4b+4)}$, then $HH_1^\s (\A)  \isom   k [ (1,A)] \oplus k [ (A^{2b+1}, A)]$. If $\l = q^{-(4b+6)}$, then $HH_1^\s (\A)  \isom k [ (A^{b+2}, A)]$.
\item For  $cd \neq 0$, $c-d \neq 0$, if $\l = q^{-4}$ then $HH_1^\s (\A)  \isom k[ (A,A) ]$.\\ 
If $\l = q^{-(2b+6)}$, then $HH_1^\s (\A) \isom k [(1,A)] \oplus k[(A^{b+2} ,A)]$.
\end{enumerate}
where for conciseness we denote by $[(x,y)]$ the class in $HH_1^\s (\A)$ of $x \otimes y \in \A^{\otimes 2}$.
\end{prop}
\begin{pf}
We have
$d_1 : \A \otimes_{\A^e} \M_1 \rto \A \otimes_{\A^e} \M_0 \isom \A$
given by 
\begin{eqnarray}
&& d_1 (a_1 \otimes e_A) = a_1 . ( A \otimes 1 - 1 \otimes A^o) = a_1 A - \s(A) a_1 = a_1 A - \mu Aa_1, \nonumber\\
&& d_1 (a_2 \otimes e_{B^\ast}) = a_2 . (B^\ast  \otimes 1 - 1 \otimes {B^\ast}^o) = a_2 B^\ast - \s(B^\ast ) a_2 = a_2 B^\ast - \l^{-1} B^\ast a_2, \nonumber\\
&& d_1 (a_3 \otimes e_B) = a_3 . (B \otimes 1 - 1 \otimes B^o) = a_3 B - \s(B) a_3 = a_3 B - \l Ba_3 \nonumber
\end{eqnarray}
($\A$ is a right $\A^e$-module via
$ a . ( t_1 \otimes {t_2}^o )= \s(t_2) a t_1$).
So $(a_1 ,a_2, a_3) \in \ker (d_1) \Leftrightarrow$
\begin{equation}
\label{a1_a2_a3}
(a_1 A - \mu A a_1) +  (a_2 B^\ast - \l^{-1} B^\ast a_2) + (a_3 B - \l B a_3)  = 0
\end{equation}
Suppose for fixed $a_3$ we have solutions 
$({a_1}', {a_2}', a_3)$, $({a_1}'', {a_2}'', {a_3})$. Then\\  $({a_1}' - {a_1}'', {a_2}' - {a_2}'', 0)$
 is a solution with $a_3 =0$, and is moreover a solution of
\begin{equation}
\label{a1_a2}
(a_1 A - \mu A a_1) +  (a_2 B^\ast - \l^{-1} B^\ast a_2)  = 0
\end{equation}
So to calculate $\ker( d_1 ) / \im( d_2 )$, we first show (Lemma \ref{a3=0_or}) that (apart from one exceptional case) for any solution 
$(a_1, a_2, a_3)$ there exists an element of $\im(d_2)$ with the same $a_3$.
This reduces the problem to solving (\ref{a1_a2}). 
Repeating this procedure, we show (Lemma \ref{a2=0_or}) that except for two special cases any solution $(a_1, a_2, 0)$ is equivalent, modulo $\im(d_2)$, to a solution $( {a_1}' , 0,0)$, which reduces the problem to solving 
\begin{equation}
\label{a1}
a_1 A - \mu A a_1  = 0
\end{equation}
Suppose for any  $a_3 = B^m A^k$ ($m \in \bZ$, $k \geq 0$) we can either find a solution $a_1 = a_1 (m,k)$, $a_2 = a_2 (m,k)$ or show that none exists. 
Let $\sS = \{ (m,k) \in \bZ \times \bN \; :$ (\ref{a1_a2_a3}) has a solution with $a_3 = B^m A^k \}$.
Then any solution of (\ref{a1_a2_a3})  is of the form 
$$a_3 = \sum_\sS \a_{m,k} B^m A^k, \; a_2 = \sum_\sS \a_{m,k} a_2 (m,k) + {a_2}', \; a_1 = \sum_\sS \a_{m,k} a_1 (m,k) + {a_1}' + {a_1}''$$
for some $\a_{m,k} \in k$, where $( {a_1}' , {a_2}')$, ${a_1}''$ are solutions of (\ref{a1_a2}), (\ref{a1}).
 
 We have
$d_2 : \A \otimes_{\A^e} \M_2 \rto \A \otimes_{\A^e} \M_1$
 given by
$$d_2 [\; 
b_1 \otimes e_A \wedge e_{B^\ast}
 + b_2 \otimes e_A \wedge e_B
 + b_3 \otimes \vartheta_S^{(1)}
+ b_4 \otimes \vartheta_T^{(1)} \;] =$$
$$[\; ( \l^{-1} B^\ast b_1 - q^2 b_1 B^\ast ) + (q^2 \l B b_2 - b_2 B) 
 -q( q^2 (b_3 A + \mu A b_3) + (d-c) b_3 )$$
\begin{equation}
\label{im_d2_1}
-q^{-1} ( b_4 A + \mu A b_4 + (d-c) b_4) \;] \otimes e_A
\end{equation}
\begin{equation}
\label{im_d2_2}
+ [ \; (b_1 A - q^2 \mu A b_1) - q^{-1} ( b_3 B + \l B b_4 ) \; ] \otimes e_{B^\ast}
\end{equation}
\begin{equation}
\label{im_d2_3}
 +[ \; (q^2 b_2 A - \mu A b_2) - q^{-1} ( \l^{-1}  B^\ast b_3 + b_4 B^\ast ) \; ] \otimes e_B
 \end{equation}

\begin{lemma}
\label{a3=0_or}
 Given 
 $(a_1, a_2, a_3) \in$ $\ker( d_1 ) / \im( d_2)$, 
we can take $a_3 =0$ unless $\l = 1$, in which case 
 the space of (equivalence classes of) solutions with $a_3 \neq 0$ is spanned (as a $k$-vector space) by $\{ a_1 = 0 = a_2, \; a_3 = B^j, \; j \geq 0 \; \}$. 
\end{lemma}
\begin{pf}
To solve (\ref{a1_a2_a3}) with $a_3 = {B^\ast}^{j+1} A^k$, take $b_3 = - q \l{B^\ast}^j A^k$, $b_1 = 0 = b_2 = b_4$ in (\ref{im_d2_3}). To solve (\ref{a1_a2_a3}) with $a_3 = B^j A^{k+1}$, take $b_2 = B^j A^k$, all other $b_i$ zero in (\ref{im_d2_3}).
This leaves the case of solving (\ref{a1_a2_a3}) with $a_3 = B^j$. Take 
$b_1 = 0$, 
$$b_2 = B^j [ q^{-2j} (1 + x^2 ) A + (d-c) (1+x) ] , \quad b_3 =  q\l ( \mu - x^2 ) B^{j+1}, \quad b_4 =  q ( \mu -x ) B^{j+1}$$
where $x = q^{2j+2}$, giving $a_3 = (x- \mu) cd B^j$.
So for $cd \neq 0$ we're done. For $cd =0$, it is clear there is no solution to (\ref{a1_a2_a3}) with $a_3 = B^j$ unless $\l =1$, in which case $a_2 = a_1 =0$ gives a solution.
\end{pf}

So we have reduced solving (\ref{a1_a2_a3}) modulo $\im( d_2 )$ to  solving (\ref{a1_a2}). In the same way, it is straightforward to show that:

\begin{lemma}
\label{a2=0_or}
 Any solution of (\ref{a1_a2_a3}) with $a_3 =0$ is equivalent, modulo $\im(d_2)$, either to a solution with $a_3 =0 = a_2$, or to one of the special cases:
\begin{enumerate}
\item $\l =1$, $\mu = \pm 1$, $a_1 =0=a_3$, $a_2 = (B^\ast)^j$, $j \geq 0$.
\item $\l =1$, $\mu = \pm 1$, $a_3 = 0$,  $a_1 =B^j [ f(2j+2) q^{-2j} A + (d-c) f(j+1)]$, 
$a_2 = ( \mu - q^{2j} ) B^{j+1}$, $j \geq 0$, which is equivalent to $a_1 = 0 = a_2$, $a_3 = B^j$.
Here $f(n) = \l^{-1} - q^{2n}$ as before. 
\end{enumerate}
\end{lemma}
 
 Finally we need to solve (\ref{a1}). For $\mu = -1$, the only solution is $a_1 =0$.
 
 \begin{lemma}
 \label{a1=0_or}
  For $\mu=1$,  $\V = \{ (a_1, 0,0) \in \ker( d_1) / \im(d_2) \}$ is spanned by:
   \begin{enumerate}
 \item If $\l \notin \{ q^{-(2b+4)} \}_{b \geq 0}$, then $(1,0,0)$ spans $\V$.
 \item $cd =0$, $\l = q^{-(2b+4)}$. Then $(A^{b+1},0,0)$ spans $\V$.
 \item $c-d =0$. For $\l = q^{-(4b+4)}$, $(A^{2b+1},0,0)$,  $(1,0,0)$  span $\V$.\\
 For $\l = q^{-(4b+6)}$, $(A^{2b+2},0,0)$ spans $\V$.
 \item $cd \neq 0$, $c-d \neq 0$. For $\l = q^{-4}$, $(A,1,1)$  spans $\V$.\\ For $\l = q^{-(2b+6)}$, $(A^{b+2},0,0)$, $( 1,0,0)$ span $\V$.
 \end{enumerate}
 \end{lemma}
 \begin{pf}
  For $\mu = 1$, the space of solutions of (\ref{a1}) is spanned by $\{ a_1 = A^j, \; j \geq 0 \}$. These solutions are not linearly independent.
 Take $b_1 = B A^j$, $b_2 =0 = b_3 = b_4$ in (\ref{im_d2_1})-(\ref{im_d2_3}), giving $a_2 = 0 = a_3$, $a_1 =   cd f( j+1) A^j + (c-d) f( j+2) A^{j+1}  - f(j+3 ) A^{j+2} $.
 Let $y_n =  f(n+1) [ A^n \otimes e_A ] \in \ker( d_1 ) / \im( d_2)$. So we have
 \begin{equation}
 \label{y_n_relation}
 y_{n+2} + (d-c) y_{n+1} -cd y_n =0 \quad \forall \; n \geq 0
 \end{equation}
 This is the same recursion relation as (\ref{x_n_recurrence_relation}).
  In addition, taking $b_1 = 0 = b_2$, $b_3 = q$, $b_4 = - q \l^{-1}$ in (\ref{im_d2_1})-(\ref{im_d2_3}), gives $a_2 = 0 = a_3$, $a_1 = 2 f(2) A + (d-c) f(1)$, hence $2 y_1 = (c-d) y_0$.
  Solving (\ref{y_n_relation}) in the same manner as for (\ref{x_n_recurrence_relation}) in the proof of Proposition \ref{HH_0}, together with this additional constraint gives the result. 
  \end{pf}

Given  $a_1 \otimes e_A + a_2 \otimes e_{B^\ast} + a_3 \otimes e_B \in \ker(d_1) / \im(d_2)$
we manufacture a twisted Hochschild 1-cycle using  (\ref{defn_f_2}). Collecting  the results of Lemmas \ref{a3=0_or}, \ref{a2=0_or}, \ref{a1=0_or} gives the description of $\ker(d_1) / \im(d_2)$ appearing in the statement of Proposition \ref{HH_1}.
This completes the proof of Proposition \ref{HH_1}.
\end{pf}

\subsection{ $HH_n^\s (\A )$, $n \geq 2$}

\begin{thm} 
\label{thm_H_2_nonzero}
\label{HH_2}
For arbitrary $c$ and $d$ (with $c+d \neq 0$), we find that:
\begin{enumerate}
\item{ For $\s = \s_\l$, $\l = q^{-(2b+2)}$, some $b \geq 0$, then $HH_2^\s (\A) \isom k$.
These automorphisms are precisely the positive powers of the modular automorphism $\smod$ (\ref{mod_aut_podles}) induced from the Haar state on quantum $SU(2)$.
}
\item{For all other $\s_\l$, $\t_\l$, $HH_2^\s (\A) =0$.}
\end{enumerate}
\end{thm}
The proof proceeds in the same manner as Proposition \ref{HH_1}, using (\ref{defn_d2}), (\ref{defn_d3}). We omit the details. 
For  $\l = q^{-(2b+2)}$,  $HH_2^\s (\A) \isom k [ \omega_2]$, where $\om_2$ is the twisted Hochschild 2-cycle:
\begin{eqnarray}
& \omega_2 = & 2[ (A^{b+1} , B, B^\ast) - ( A^{b+1} , B^\ast ,B) + 2( A^b B, B^\ast , A) - 2 q^{-2} (A^b B, A, B^\ast )] \nonumber\\
&& + 2(q^4 -1) ( A^{b+1} , A,A) + (1 - q^{-2} ) cd (c-d) ( A^b ,1,1) \nonumber\\
&& \label{omega_2}
+(c-d)[ (A^b , B^\ast , B) - q^{-2} ( A^b , B, B^\ast ) + (1- q^2 ) ( A^b , A,A) ] 
\end{eqnarray}

Finally, all the higher twisted Hochschild homology groups vanish:

\begin{prop}
We have $HH_n^\s (\A) =0$ for all $n \geq 3$ for any  $\s$.
\end{prop}

We prove this in the case $n=3$:

\begin{thm}
\label{HH3=0}
 $HH_3^\s (\A )=0$ for any automorphism $\s$.
\end{thm}
\begin{pf}
We have 
$$d_3 [ a_1 \otimes (e_A \wedge \vartheta_S^{(1)}) + a_2 \otimes  (e_A \wedge \vartheta_T^{(1)}) + a_3 \otimes (e_{B^\ast} \wedge \vartheta_S^{(1)}) + a_4 \otimes  (e_B \wedge \vartheta_T^{(1)}) ] $$
\begin{equation}
\label{im_d3_1}
= [ q^{-3} a_1 B + q^{-1} \l B a_2 - q^{-1} ( a_3 A + q^2 \mu A a_3 - (c-d) a_3)  ] \otimes (e_A \wedge e_{B^\ast})
\end{equation}
\begin{equation}
\label{im_d3_2}
+ [ q^{-3} \l^{-1} B^\ast a_1 + q^{-1} a_2 B^\ast - q^{-1} ( q^2 a_4 A + \mu A a_4 - (c-d) a_4 ) ] \otimes( e_A \wedge e_B)
\end{equation}
\begin{equation}
\label{im_d3_3}
+ [ ( a_1 A - \mu A a_1 ) + a_3 B^\ast - \l B a_4  ] \otimes \vartheta_S^{(1)}
\end{equation}
\begin{equation}
\label{im_d3_4}
+ [ (a_2 A - \mu A a_2 ) - \l^{-1} B^\ast a_3 + a_4 B ] \otimes \vartheta_T^{(1)}
\end{equation}
and $d_4 [ 
b_1 \otimes (e_A \wedge e_{B^\ast} \wedge \vartheta_S^{(1)})
+ b_2 \otimes (e_A \wedge e_B \wedge \vartheta_T^{(1)}) 
+ b_3 \otimes \vartheta_S^{(2)}
+ b_4 \otimes \vartheta_T^{(2)} ] = $
\begin{eqnarray}
\label{im_d4_2}
&&=[ - q^2 b_1 B^\ast + q^2 \l B b_2 - q( q^2 b_3 A + \mu q^2 A b_3 - (c-d) b_3 )] \otimes (e_A \wedge \vartheta_S^{(1)})\nonumber\\
&&+[ \l^{-1} B^\ast b_1 - b_2 B - q^{-1} ( b_4 A + \mu A b_4 - (c-d) b_4 ) ] \otimes (e_A \wedge \vartheta_T^{(1)})\nonumber\\ 
\label{im_d4_3}
&&+[ ( b_1 A - q^2 \mu A b_1 ) - q^{-1} b_3 B - q^{-1} \l B b_4 ] \otimes ( e_{B^\ast} \wedge \vartheta_S^{(1)})\nonumber\\ 
\label{im_d4_4}
&&+[ ( q^2 b_2 A - \mu A b_2 ) - q^{-1} \l^{-1} B^\ast b_3 - q^{-1} b_4 B^\ast  ] \otimes ( e_B \wedge \vartheta_T^{(1)})
\end{eqnarray}

Finding $\ker( d_3)$ corresponds to finding all solutions $(a_1, a_2, a_3, a_4)$ to the four equations (\ref{im_d3_1})-(\ref{im_d3_4}).
Our strategy is the same as for Proposition \ref{HH_1}. Suppose for fixed $a_4$ we  find solutions $(a_1, a_2, a_3, a_4)$, $({a_1}', {a_2}', {a_3}', a_4)$. Then  $(a_1 - {a_1}', a_2 - {a_2}', a_3 - {a_3}', 0)$ is a solution with $a_4 =0$. 
So to calculate $\ker( d_3 ) / \im( d_4 )$, we first show (Lemma \ref{lemma_a4=0}) that for any solution 
$(a_1, a_2, a_3, a_4)$ there exists an element of $\im(d_4)$ with the same $a_4$.
So  we need only look for solutions with $a_4 =0$.

We  repeat this procedure for $a_3$ (Lemma \ref{lemma_a3=0}), showing that $\ker(d_3) / \im(d_4)$ is spanned by (equivalence classes of) solutions with $a_3 =0 = a_4$. 
Finally we show (Lemmas \ref{lemma_a1=a2=0_mu=-1}, \ref{lemma_a1=a2=0_mu=1}) that any such solution belongs to $\im(d_4)$. 

\begin{lemma}
\label{lemma_a4=0} Any solution $(a_1, a_2, a_3, a_4)$ of (\ref{im_d3_1})-(\ref{im_d3_4}) is equivalent, modulo $\im(d_4)$, to a solution with  $a_4 =0$.
\end{lemma}
\begin{pf}
We start by solving for given $a_4$. It is enough just to consider monomials. 
For $a_4 = ( B^\ast )^{j+1} A^k$, 
take $b_3 = - q \l( B^\ast )^{j} A^k$, 
$b_2 = 0 = b_4$  in (\ref{im_d4_4}). 
To solve for $a_4 = B^j A^{k+1}$, take 
$b_2 = ( q^2 - \mu q^{-2j} )^{-1} B^j A^k$, $b_3 = 0 = b_4$.
 Then $q^2 b_2 A - \mu A b_2 =  B^j A^{k+1}$.
So we are left with only the case $a_4 = B^j$. 
Take $b_2 = ( q^2 - \mu q^{-2j} )^{-1} B^j [ \a_0 + \a_1 A]$, $b_3 = -q B^{j+1}$, $b_4 =0$.
Then $- q^{-1} b_4 B^\ast - q^{-1} \l^{-1} B^\ast b_3 = B^j [ cd + q^2 (c-d) A + q^4 A^2 ]$
and $q^2 b_2 A - \mu A b_2 = B^j [ \a_0  A + \a_1 A^2 ]$. Taking
$\a_0 = q^2 (d-c)$,  $\a_1 = - q^4$,
we see that provided $cd \neq 0$, we can find solutions with $a_4 = B^j$ for any $j \geq 0$.
If $cd = 0$, then we see from (\ref{im_d3_3}), (\ref{im_d3_4}) that $a_4 = B^j$ cannot be in $\ker(d_3)$.
\end{pf}

In the same way, it is straightforward to show that:

\begin{lemma}
\label{lemma_a3=0}
Any solution $(a_1, a_2, a_3, a_4)$ of (\ref{im_d3_1})-(\ref{im_d3_4}) with $a_4 =0$ is equivalent, modulo $\im(d_4)$, to a solution with   $a_3 = a_4 =0$.
\end{lemma}

So we need only consider $a_1$, $a_2 \neq 0$. From (\ref{im_d3_3}), (\ref{im_d3_4}), we have
\begin{equation}
\label{a1_a2_first_condition}
a_1 A = \mu A a_1, \quad a_2 A = \mu A a_2
\end{equation}
\begin{lemma}
\label{lemma_a1=a2=0_mu=-1}
 For $\mu = -1$, the only solution to (\ref{a1_a2_first_condition}) is $a_1 = 0 = a_2$. 
\end{lemma}
Hence for $\mu = -1$, $\ker( d_3 ) = \im( d_4)$, thus proving Theorem \ref{HH3=0} in this case.

For $\mu =1$, (\ref{im_d3_1}), (\ref{im_d3_2}) give 
$a_1 B + q^2 \l B a_2 =0$, 
$B^\ast a_1 + q^2 \l a_2 B^\ast =0$
(it is straightforward to show that these two conditions are equivalent).
So for $\mu = 1$, $\ker(d_3) / \im( d_4)$ is spanned by (the equivalence classes of) the solutions
\begin{equation}
\label{a1_a2_nonzero}
a_1 = - \l q^{2j+2} A^j, \quad a_2 = A^j, \quad a_3 = 0 = a_4 \quad (j \geq 0)
\end{equation}
\begin{lemma}
\label{lemma_a1=a2=0_mu=1}
 The solutions (\ref{a1_a2_nonzero}) all belong to $\im(d_4)$.
\end{lemma}
\begin{pf}
In the case $cd \neq 0$, $c-d \neq 0$, taking
$$b_1 = 4 \a_1 B A^j, \quad b_2 = 4 q^{2j} ( \a_1 \l^{-1} - \g) B^\ast A^j, \quad \g = 4 (c+d)^{-2}$$
$$b_3 =  \l \g q^{2j+1} A^j [ 2 q^2 A - (c-d) ], \quad b_4 = -  \g q A^j [ 2A - (c-d)]$$
some $\a_1 \neq \l \g$, gives  (\ref{a1_a2_nonzero}).
The other two  cases ($cd =0$, $c=d$) are similar. 
\end{pf}
This completes the proof of Theorem \ref{HH3=0}. 
\end{pf}

\section{Twisted cyclic homology of the Podle\'s spheres}
\label{section:twisted_cyclic_homology}

For an  algebra $\A$ and automorphism $\s$, twisted cyclic homology $HC_*^\s (\A)$ arises as in \cite{loday} from the 
cyclic module $\Asn$, with objects ${\{ \Asnn \}}_{n \geq 0}$ (\ref{C_n^sigma}) defined  by 
$\Asnn = \A^{\otimes (n+1)} / (\id - \s^{\otimes (n+1)})$. 
The face, degeneracy and cyclic operators were given explicitly in \cite{hk}.
Twisted cyclic homology $HC_*^\s (\A)$ is the total homology of Connes' mixed $(b,B)$-bicomplex
corresponding to the cyclic module $\Asn$:
\begin{equation}
\label{mixed_complex}
     	\begin{CD}
	@ V{b_4} VV @ V{b_3} VV @ V{b_2} VV @ V{b_1} VV @ . @ . @ .\\
	{\Asnthree} @ <{B_2}<< {\Asntwo} @ <{B_1} << {\Asnone} @ <{B_0} << {\Asnzero} @ . @ . @ . @ .\\
 	@ V{b_3} VV @ V{b_2} VV @ V{b_1} VV @ . @ . @ . @ .\\
	{\Asntwo} @ <{B_1} << {\Asnone} @ <{B_0} << {\Asnzero} @ . @ . @ . @ . @ .\\
 	@ V{b_2} VV @ V{b_1} VV @ . @ . @ . @ . @ .\\
{\Asnone} @ <{B_0} << {\Asnzero} @ . @ . @ . @ . @ . \\
 	@ V{b_1} VV @ . @ .  @ . @ . @ . @ .\\
	{\Asnzero} @ . @ . @ . @ . @ . @ .\\
    	\end{CD}
\end{equation}
The  maps $b_n$ coincide with the twisted Hochschild boundary maps $b_\s$ (\ref{b_sigma}).
We will drop the suffices and write $b_n$, $b_\s$ as  $b$.
In lowest degrees, the maps $B_n$
are:
$$B_0 [a_0] = [(1,a_0)] + [ (\s(a_0), 1)] = [(1, a_0 )] + [(a_0,1)],$$
$$B_1 [(a_0, a_1)] =  [(1,a_0, a_1)] - [(\s(a_1), 1, a_0 )] - [(1, \s(a_1), a_0 )] + [(a_0, 1, a_1)]$$
For any $a \in \A$, $b(a,1,1) = (a,1)$, so $[(a,1)] = 0$ in $HH_1^\s (\A)$.
So the induced map $B_0 : HH_0^\s (\A) \rto HH_1^\s (\A)$ satisfies
$B_0 [a] = [(1,a)]$.
For $t \in \A$, with $\s(t) = \a t$, some $\a \in k$, then
\begin{equation}
\label{first_homology_relation}
b ( \Sigma_{j=0}^m \; \a^j (t^j , t^{m-j}, t) - (t^{m+1}, 1,1) )=
(\Sigma_{j=0}^m \; \a^j) \; (t^m , t) - (1, t^{m+1})
\end{equation}
If $\a = 1$, then 
$B_0 [ t^{m+1}] = [(1,t^{m+1})] = (m+1) [(t^m ,t)] \in HH_1^\s (\A)$.

Taking $\A = \Acd$, we calculate total homology of the mixed complex (\ref{mixed_complex}) via a spectral sequence.
The first step (vertical homology of the columns) gives: 
\begin{equation}
\label{first_page}
     	\begin{CD}
	@ VVV @ VVV @ VVV @ VVV @ . @ . @ .\\
	{0} @ <{0}<< {HH_2^\s (\A)} @ <{B_1} << {HH_1^\s (\A)} @ <{B_0} << {HH_0^\s (\A)} @ . @ . @ . @ .\\
 	@ VVV @ VVV @ VVV @ . @ . @ . @ .\\
	{HH_2^\s (\A)} @ <{B_1} << {HH_1^\s (\A)} @ <{B_0} << {HH_0^\s (\A)} @ . @ . @ . @ . @ .\\
 	@ VVV @ VVV @ . @ . @ . @ . @ .\\
{HH_1^\s (\A)} @ <{B_0} << {HH_0^\s (\A)} @ . @ . @ . @ . @ . \\
 	@ VVV @ . @ .  @ . @ . @ . @ .\\
	{HH_0^\s (\A)} @ . @ . @ . @ . @ . @ .\\
    	\end{CD}
\end{equation}
since for every $\s$ we have $HH_n^\s (\A) =0$ for $n \geq 3$.
We find that:

\begin{prop}
\label{generic_lambda}
For $\l \notin q^{-2 \bN}$, $\s = \s_\l$,
$HC_{2n}^\s (\A) = k [1] \oplus k [A]$,\\
$HC_{2n+1}^\s (\A) =0$. For $\s = \t_\l$ with $\l \neq 1$, $HC_n^\s (\A) =0$ for $n \geq 1$. 
\end{prop}
\begin{pf} In both cases $HH_0^\s (\A) = k [1] \oplus k [A]$ (with $[A] =0$ for $\s = \t_\l$), 
 $HH_n^\s (\A) = 0$ for $n \geq 1$, 
   (\ref{first_page}) stabilizes immediately, and the result follows. 
\end{pf}
\begin{prop}
For $\l =1$, $\mu = \pm 1$, then just as in \cite{mnw} we have 
$$HC_0^\s (\A) =  k [1] \oplus k [A] \oplus \; ( \Sigma_{m>0}^\oplus \; k [B^m] ) \; \oplus \; ( \Sigma_{m>0}^\oplus \; k [{B^\ast}^m] )$$
$HC_{2n+1}^\s (\A) =0, \quad
HC_{2n+2}^\s (\A) = k [1] \oplus k [A], \quad$
with $[A] =0$ for $\mu = -1$.
\end{prop}
\begin{pf} We have
$B_0 [1] =0$, $B_0 [A] = [(1,A)] =0$, 
while $B_0 [ B^{m+1} ] =  [(1,B^{m+1}) ]$ $= (m+1) [ (B^m , B)]$ by (\ref{first_homology_relation}), and in the same way 
$B_0 [ {B^\ast}^{m+1} ] = (m+1) [ ({B^\ast}^m , B^\ast )]$.
So $\ker (B_0) = k[1] \oplus k[A]$, and $HH_1^\s (\A) = \im (B_0)$. 
Hence the spectral sequence stabilizes at the second page with all further maps being zero. 
\end{pf}

\begin{prop} For $\s = \s_\l$, $\l = q^{-(2b+2)}$, then $HC_{2n+1}^\s (\A) = 0$, and:
\begin{enumerate}
\item $\l = q^{-2}$. 
 $HC_{2n+2}^\s (\A) = k[1] \oplus k[ \om_2]$.
\item $\l = q^{-4}$. For
 $c =d$, $HC_{2n+2}^\s (\A) = k[ \om_2]$,
  else $HC_{2n+2}^\s (\A) = k[1] \oplus k[ \om_2]$.
\item \label{we_prove_this_one}
$\l = q^{-(4b+6)}$.
For $cd =0$ or $c=d$,  
 $HC_{2n+2}^\s (\A) = k[1] \oplus k[ \om_2]$, otherwise
 $\quad HC_{2n+2}^\s (\A) = k[ \om_2]$.
\item $\l = q^{-(4b+8)}$. 
For $cd =0$, $HC_{2n+2}^\s (\A) = k[1] \oplus k[\om_2]$, 
otherwise\\ $HC_{2n+2}^\s (\A) = k[ \om_2]$.
\end{enumerate}
\end{prop}
\begin{pf} We prove case \ref{we_prove_this_one}, the others are completely analogous.
For $cd =0$ or $c=d$, $HH_0^\s (\A) = k[1] \oplus k[A^{2b+3}]$, and $HH_1^\s (\A) = k[(A^{b+2}, A]$.
We have $B_0 [1] = [(1,1)] =0$, $B_0 [A^{2b+3}] = [(1, A^{2b+3})] = (2b+3) [(A^{b+2}, A)]$.
So $\ker( B_0 ) = k[1]$, $\im( B_0) = HH_1^\s (\A)$. Then
the spectral sequence
 (\ref{first_page}) stabilizes at page two: 
	\[
     	\begin{CD}
	@ VVV @ VVV @ VVV @ VVV @ . @ . @ .\\
	{0} @ <<< {k[\omega_2]} @ <<< {0} @ <<< {k[1]} @ . @ . @ . @ .\\
 	@ VVV @ VVV @ VVV @ . @ . @ . @ .\\
	{k [ \omega_2]} @ <<< {0} @ <<< {k[1]} @ . @ . @ . @ . @ .\\
 	@ VVV @ VVV @ . @ . @ . @ . @ .\\
{0} @ <<< {k[1]} @ . @ . @ . @ . @ . \\
 	@ VVV @ . @ .  @ . @ . @ . @ .\\
	{HH_0^\s (\A)} @ . @ . @ . @ . @ . @ .\\
    	\end{CD}
	\]
with all further maps being zero. 
For $cd \neq 0$ and $c \neq d$, then $HH_0^\s (\A) = k[1] \oplus k[A^{2b+3}] =  k[A] \oplus k[A^{2b+3}]$, and $HH_1^\s (\A) = k[(1,A)] \oplus k[(A^{b+2}, A]$. Then $B_0 [A] = [(1,A)]$, hence $\ker( B_0) =0$, $\im( B_0) = HH_1^\s (\A)$.
\end{pf}

\section{The standard Podle$\acute{\mathrm{s}}$ quantum sphere}
\label{standard_quantum_sphere}

We specialize our results to the standard quantum sphere $\A(S^2_q)$, which as described in Section \ref{section:podles} naturally embeds as a *-subalgebra of $\A(SU_q(2))$.  
We recall that Schm\"udgen and Wagner \cite{sw2} defined a twisted cyclic 2-cocycle $\tau$  over $\A(S^2_q)$ as follows.
For $a_0$, $a_1$, $a_2 \in \Asq$, define
\begin{equation}
\label{defn_tau}
\tau ( a_0 , a_1 ,a_2 ) = h ( a_0 [  (a_1 \tl F) (a_2 \tl E) - q^2 (a_1 \tl E)  ( a_2 \tl F) ] )
\end{equation}
where $\tl$ is the right action of $U_q ( su(2) )$  (\ref{right_action}). 
As shown in \cite{sw2}, the mappings $\A( S^2_q) \rto \A( SU_q (2))$
given by $x \mapsto  x \tl E$, $x \mapsto x \tl F$ are derivations.
Here $h$ denotes the Haar state on $\A(SU_q (2))$, 
 which restricts to $\Asq$ as
$$ h ( A^r B^s) = 0 = h ( A^r (B^* )^s) \quad s >0,\quad
 h( A^r ) = (1- q^2 ) (1 - q^{2r+2} )^{-1}$$
Schm\"udgen and Wagner proved:

\begin{prop}
\cite{sw2}, Theorem 4.5.
$\t$ is a nontrivial $\s$-twisted cyclic 2-cocycle on $\Asq$, with $\s$ the automorphism given by $\s(x) = K^{-2} \triangleright x$.  Further, $\t$ is $ U_q (su(2))$-invariant  and coincides with the volume form of the distinguished covariant 2-dimensional first order differential calculus found by Podle\'s  \cite{pod92}.
\end{prop}
Schm\"{u}dgen and Wagner also constructed a $U_q (su(2))$-equivariant Dirac operator, unitarily equivalent to those previously found by Bibikov and Kulish \cite{bk} and Dabrowski and Sitarz \cite{ds}, which they used to give a representation of the Podle\'s calculus by  bounded commutators.

Explicitly,  $\s (B) = q^2 B$, $\s (B^* ) = q^{-2} B^*$.
So in (\ref{defn_sigma}),  $\l = q^2$. From Proposition \ref{HH_1} and Theorem \ref{HH_2} we have $HH_n^\s (\A) =0$ for $n \geq 1$ for this $\s$, i.e. this twisted cocycle does not correspond to the ``no dimension drop" case. 
By Proposition \ref{generic_lambda}, we have $HC_{2n}^\s (\A) = \bC [1] \oplus \bC [A]$, $HC_{2n+1}^\s (\A) =0$ for all $n \geq 0$. 
The $\s$-twisted cyclic 0-cocycles $\t_0$, $h_A$ dual to $[1]$, $[A]$ are defined on Poincar\'e-Birkhoff-Witt monomials $x$  (\ref{pbw_basis}) by
$\t_0 (1) =1$, $\t_0 (x) =0$ for  $x \neq 1$, and 
$$ h_A ( A^r B^s ) = 0 = h_A ( A^r (B^* )^s  ) \quad s > 0$$
$$h_A (1) =0, \quad h_A (A^{r+1}) = (1 - q^4 )(1 - q^{2r+4} )^{-1} $$
The Haar state $h$ (restricted to $\A(S^2_q)$) is given by $h = \t_0 +  (1+ q^2)^{-1} h_A$.
By cohomology calculations completely dual to our previous homology calculations, we have $HC^{2n}_\s (\A) \isom \bC [ S^n \t_0 ] \oplus \bC [ S^n h_A]$, $HC^{2n+1}_\s (\A) =0$, where $S$ is Connes' periodicity operator.
We can now identify the class of $\t$ in $HC^2_\s (\A)$:

\begin{thm} 
We have $[\t] = \b [ Sh_A ] \in HC^2_\s (\A)$, for some nonzero $\b$.
\end{thm}
\begin{pf}
We have $HC^2_\s (\A) \isom \bC^2$, generated by $[ S \t_0 ]$, $[S h_A ]$, where $S \phi ( a_0 , a_1 , a_2 ) = \phi ( a_0 a_1 a_2 )$ for any $\phi \in HC^0_\s (\A)$.
Recall from \cite{sw2} the element
$$\eta = (B^* ,A,B) + q^2 (B, B^* ,A) + q^2 (A,B,B^* ) - q^{-2} (B^* ,B,A) $$
$$- q^{-2} (A, B^* ,B) - (B,A, B^* ) + (q^6 - q^{-2}) (A,A,A)$$
Now, $\t (\eta) = -1$, and it was shown in  \cite{sw2} that $[\t]$ is nontrivial in $HC^2_\s (\A)$. 
So there are scalars $\a$, $\b$, not both zero,  such that $[ \t ] = \a [ S \t_0 ] + \b [ S h_A ]$.
Now, $\t (1,1,1)=0 = Sh_A (1,1,1)$, whereas $S \t_0 (1,1,1) = \t_0 (1) =1$. Hence $\a =0$.
Since $S \eta = (q^4 - q^{-2}) A^2$, we have $S h_A ( \eta) = h_A ( S \eta) = (q^4 - q^{-2}) h_A (A^2) = q^2 - q^{-2}$.
 If $\eta$ was a twisted 2-cycle we could  deduce that $\b = ( q^{-2} - q^2)^{-1}$.
Since  $b_\s (\eta) = 2( q^4 - q^{-2}) (A,A) \neq 0$, this need not  hold. 
We could calculate $\b$ by finding $[ {\bf a} ]  \in HC_2^\s (\A)$ such that $[ S {\bf a} ] = [A] \in HC_0^\s ( \A)$
(note that $(1 - q^{2s+4} )[ A^{s+1}] = ( 1- q^4 ) [A]$ for $s \geq 0$).
 Then 
 $\t ( {\bf a} ) = \b Sh_A ( {\bf a}) = \b h_A ( S{\bf a} ) = \b h_A (A) =  \b$.
However finding such an ${\bf a}$ explicitly has not been possible. 
\end{pf}
\section{Acknowledgements}

I am grateful for  the support of the EU Quantum Spaces 
network (INP-RTN-002) and of the EPSRC via a Postdoctoral Fellowship.
I also thank Ulrich Kr{\"a}hmer and the referee for their very useful comments.

\end{document}